\documentclass[a4paper,11pt]{article}   % MikTeX
\usepackage{lineno,hyperref}
\usepackage[utf8]{inputenc}
\usepackage[english]{babel}
\usepackage{ifthen}
\usepackage{amssymb,amsmath}
\usepackage[usenames]{color}
\date{}

 \newif\ifNoRemark
    \def\addtheorem#1#2#3#4{ % \usepackage{ifthen} needed
    \ifthenelse{\expandafter\isundefined\csname the#2\endcsname}{\newcounter{#2}}{}
    \newenvironment{#1}[1][\global\NoRemarktrue]% No Remark by default
     {\par\addvspace{2mm}\noindent
       \refstepcounter{#2}{\bf #3~\csname the#2\endcsname
      \vphantom{##1}\ifNoRemark.\ \else\ (##1).\fi}\begingroup #4}%
     {\endgroup\par\addvspace{1mm}\global\NoRemarkfalse}
    \expandafter\newcommand\csname b#1\endcsname{\begin{#1}}
    \expandafter\newcommand\csname e#1\endcsname{\end{#1}}
    }

\addtheorem{theorem}{thrm}{Theorem}{\sl}
\addtheorem{lemma}{lmm}{Lemma}{\sl}
\addtheorem{proposition}{prpstn}{Proposition}{\sl}
\addtheorem{corollary}{crllr}{Corollary}{\sl}
\addtheorem{problem}{problem}{Problem}{\sl}
\addtheorem{remark}{}{Remark}{}

\begin{document}

\title{Embedding in MDS codes and Latin cubes\footnote{\,The study
was carried out within the framework of the state contract of the
Sobolev Institute of Mathematics (project no. 0314-2019-0017).}  }

\author{ Vladimir N. Potapov\\
Sobolev Institute of Mathematics, Novosibirsk, Russia \\
vpotapov@math.nsc.ru\\}

\maketitle

\begin{abstract}
An embedding of a code is a mapping that preserves distances between
codewords. We prove that any code with code distance $\rho$ and
length $d$ can be embedded into an MDS code with the same code
distance and  length but under a larger alphabet. As a corollary we
obtain embeddings of systems of partial mutually orthogonal Latin
cubes and $n$-ary quasigroups.


\end{abstract}


\section{Introduction}

Let $Q_q=\{0,\dots,q-1\}$. Consider the Cartesian product  $Q_q^d$
as the metric space with the Hamming distance $\rho$.   A subset $C$
of $Q_q^d$ is called a code with distance $t+1$ if $\min
\rho(x,y)=t+1$ for any $x,y\in C$, $x\neq y$. A subset $C$ of
$Q_q^d$ is called an $MDS(t,d,q)$ code (of order $q$, code distance
$t+1$ and length $d$) if $|C\cap\Gamma|=1$ for each $t$-dimensional
axis-aligned plane $\Gamma$. Ethier and Mullen \cite{EMullen} proved
that $MDS(t,t+s,q)$ codes are equivalent to a set of $t$ mutually
orthogonal  $s$-dimensional Latin cubes of order $q$. If $s=2$ then
an $MDS(t,t+2,q)$ code is equivalent to a system of MOLS; if $t=1$
then an $MDS(1,1+s,q)$ code is equivalent to an $s$-dimensional
Latin cube or the Cayley table of an $s$-ary quasigroup of order
$q$. The main result of the present paper is the following

\begin{theorem}\label{thEmbed1}
Let $C\subset Q_q^d$ be a code with distance $t+1$. Then there exist
$q', q'\geq q$ and  an $MDS(t,d,q')$ code $M$ such that $C\subseteq
M\subset Q_{q'}^d$.
\end{theorem}

It means that $Q_q\subseteq Q_{q'}$ and $Q_q^d$ is a subset of
$Q_{q'}^{d}$ where $q\leq q'$. Consequently, if $C\subset Q_q^d$
then $C\subset Q_{q'}^{d}$.

This theorem generalizes theorems on embeddings    of partial Latin
squares \cite{Evans}, $d$-dimensional Latin cubes \cite{Cruse} and
systems of mutually orthogonal partial Latin squares \cite{Linder},
\cite{DonY}. These embeddings preserves dimension $d$ and increases
order $q$. In Theorem \ref{thEmbed1} the order (size of the
alphabet) of a code increases exponentially. In a construction
 for embedding partial MOLS  from \cite{DonY} the order of code increases
 polynomially. In Theorem \ref{thEm3} we propose a construction with
 a polynomial grow of order for embeddings of partial $d$-dimensional Latin
 cubes. Note than in the special case of Latin squares our order
 $q'=q^2$ is worse than  $q'=2q$ in classical paper by Evans
 \cite{Evans}.

 Another embedding that preserves the order is considered
in \cite{KrotovSot} by Krotov and Sotnikova. More precisely, they
proved that a code $C\subset Q_q^d$ with distance $3$ can be
embedded in a perfect $1$-error-correcting code in $Q_q^{d'}$, $d'=
\frac{q^d-1}{q-1}$, where $q$ is a prime power.

In \cite{DonGY} and \cite{DonGY1} Donovan, Grannell and Yazici
consider a more complicated problem. In  terms of the present paper,
they construct  MDS codes with partially predetermined projections.

\section{Preliminaries}

The proofs of propositions from this section can be found in
\cite{DK} and \cite{MS}. Notation and propositions are also examined
in detail in \cite{Pot18} and \cite{Pot20}.

As defined above a subset $C$ of $Q_q^d$ is called an $MDS(t,d,q)$
code (of order $q$, code distance $t+1$ and length $d$) if
$|C\cap\Gamma|=1$ for each $t$-dimensional axis-aligned plane
$\Gamma$. It means that $Q_q^d$ is an MDS code with code distance
$1$. Note that the alphabetical list is not important in the
definition of MDS codes. Moreover, we can use different alphabets in
different coordinates if it is convenient.

\begin{proposition}[Singleton bound]\label{proNM5}
A subset $M\subset Q_q^d$ with code distance $t+1$ is an {\rm MDS}
code if and only if $|M|=q^{d-t}$.
\end{proposition}

Let $C\subset Q_q^d$. The set
$$\widehat{C^i_a}=\{(x_1,\dots,x_{i-1},x_{i+1},\dots x_{d}) :   (x_1,\dots,x_{i-1},a,x_{i+1},\dots,x_{d})\in C\}.$$
is called a retract of $C$. Note that $\widehat{C^i_a}\subset
Q_q^{d-1}$.

% We  use notation
%$(C)^{x_{i_1},\dots,x_{i_k}}_{a_1,\dots,a_k}$ to fix several
%variables.

An embedded retract of  $C$ is the set
$$C^i_a=\{(x_1,\dots,x_{i-1},x_i,x_{i+1},\dots x_{d})\in C : x_i=a\}$$
Note that $C^i_a\subset Q_q^{d}$.

 A projection of  $C$ along the  $i$th direction is  the set
$$C_i=\{(x_1,\dots,x_{i-1},x_{i+1},\dots x_{d}) : \exists\,a\in Q_q\
 (x_1,\dots,x_{i-1},a,x_{i+1},\dots,x_{d})\in C\}.$$ Note that
$C_i\subset Q_q^{d-1}$.

\begin{proposition}\label{progs1}
\begin{itemize}
\item Every projection of an MDS code with distance $t+1$, $t\geq 1$, is
an MDS code with code distance $t$.
\item Every retract of an MDS code with distance $t<d$ is an MDS code
with the same code distance.
\end{itemize}
\end{proposition}

\begin{proposition}\label{proNMq}
Let $M\subset Q_q^d$ be an $MDS(1,d,q)$ code. Then the set
$$\{(x_1,\dots,x_{d-1},x_{d}+a({\rm mod}\,q), a) : a\in Q_q, (x_1,\dots,x_{d})\in M\}$$ is
an $MDS(1,d,q)$ code.
\end{proposition}

Therefore, any $MDS(1,d,q)$ code can be included as a retract in an
$MDS(1,d',q)$, where $d'>d$. An analogous property for MDS codes
with
 code distance greater than $2$ is not true.

A subset $T$ of an MDS code $M\subset Q_q^d$ is called an { MDS
subcode} of the code if $T$ is an MDS code  in $A_1\times\dots\times
A_d$ with the same code distance as $M$ and $T=M\cap
(A_1\times\dots\times A_d)$ where $A_i\subset Q_q$, $i\in
\{1,\dots,d\}$. Obviously $|A_1|=\dots=|A_d|=q'$ and $q'$ is the
order of the subcode $T$. Note that the definition of a Latin
subsquare is analogous.

\begin{proposition}\label{proNMa}%{\rm(\cite{Pot18})}
Suppose that $C$ is an MDS code, $C_1$ is an MDS subcode of $C$ in a
subcube $A_1\times\dots\times A_d$, $C_2$ is an MDS code in a
subcube $A_1\times\dots\times A_d$ with the same distance as $C_1$.
Then if we switch $C_1$ by $C_2$  then we  obtain the MDS code $C'$
with the same parameters as $C$.
%There exist two codes which are different only by the subcode.
\end{proposition}

This exchanging of subcodes is called a switching of $C$.

\begin{corollary}\label{proNMf}
Suppose that $C$ is an MDS code, $C_1$ is an MDS  subcode of $C$ in
a subcube $A_1\times\dots\times A_d$ and $u\in A_1\times\dots\times
A_d$. Then there exists an MDS code $C'=(C\setminus C_1)\cup C_2$
such that $u\in C'$.
\end{corollary}

\begin{proposition}\label{proNMh}
If $C_1$ and $B_1$ are disjoint MDS subcodes of $C$ and $C_2$ is an
MDS code in the same subcube as $C_1$ then $B_1$ is an MDS subcode
of $(C\setminus C_1)\cup C_2$.
\end{proposition}

The set $Q_{q_1q_2}$ can be considered as the Cartesian product
$Q_{q_1}\times Q_{q_2}$. Consequently, we can identify
$Q^d_{q_1}\times Q^d_{q_2}$ and the hypercube $Q^d_{q_1q_2}$. Thus,
if $C_1\subset Q^d_{q_1}$ and $C_2\subset Q^d_{q_2}$ then
$$C_1\times C_2=\{((x_1,y_1),(x_2,y_2),\dots,(x_d,y_d)) :
(x_1,\dots,x_d)\in C_1, (y_1,\dots,y_d)\in C_2\}\subset
Q^d_{q_1q_2}.$$

\begin{theorem}[McNeish]\label{proNMb1}
Suppose  $M_1$ is an  $MDS(t,d,q_1)$ code and $M_2$ is an
$MDS(t,d,q_2)$ code. Then $M_1\times M_2$ is  an  $MDS(t,d,q_1q_2)$
code.
\end{theorem}

It is  well known the following generalization of McNeish's theorem.

\begin{proposition}\label{Cartprod1}
Let $M$ be an  $MDS(t,d,q_1)$ code  and  let $U[x]$ be an
 $MDS(t,d,q_2)$ code for each $x\in M$. Then the set $\bigcup\limits_{x\in M} x
\times U[x]$ is an $MDS(t,d,q_1q_2)$ code.
\end{proposition}

Let $q$ be a prime power and let $Q_q=GF(q)$.  A linear
$k$-dimensional subspace $C\subset Q_q^d$ with  distance $t$ is
called $[d,k,t]_q$ code over $GF(q)$. By Proposition \ref{proNM5},
we see that any $[d,d-t,t+1]_q$ code over $GF(q)$ is an $MDS(t,d,q)$
code.

For each linear code $C$ there exists a check matrix $A_C$ such that
$$ C=\{x\in (GF(q))^d : A_Cx=\bar 0\}, $$
where $A_C$ is a matrix of size $(d-k)\times d$.

\begin{proposition}\label{proNMe}
A linear code $C\subset Q_q^d$ is an MDS code if and only if  all
minors of $A_C$ of order $d-k$ is nonzero.
\end{proposition}

In this case $C$ is $MDS(d-k+1,d,q)$ code. It is easy to see that we
can choose check matrix of special type $A_C=(I|A')$ where $I$ is
indentity matrix of size $t\times t$, $t=d-k$.

\begin{corollary}\label{proNMt}
Let $A_C$ be the check matrix  of an MDS code with size $(d-k)\times
d$. Then every nontrivial linear combination of rows of $A_C$
contains greater than $k+1$ nonzero elements.
\end{corollary}

\begin{proposition}\label{proNMb}
Let $q$ be a prime power. Then
 for each  integer $d\leq q+1$  and $\rho$,
$3\leq \rho < d$, there exists a linear (over $GF(q)$) {\rm MDS}
code $C\subset Q_q^d$ with the code distance $\rho$.
\end{proposition}

 We will consider  the Cartesian
product of linear codes. Let $V= (GF(q))^m$. We consider elements of
$V^d$ as matrices of size $d\times m$ over $GF(q)$. Let $A$ be the
check matrix of an MDS code $C$ over $GF(q)$. Then
\begin{equation}\label{eqEmb1}
 M=\{Y\in V^d:
AY=\bar 0\}
\end{equation}
 is an MDS code by the McNeish's theorem.
Moreover, for a linear subspace $W\subset V^d$ the set $
M|_{W}=\{Y\in W^d: AY=\bar 0\} $ is an MDS code and $ M|_{W}=M\cap
W^d$ is an MDS subcode of $M$.

\begin{proposition}\label{proNMg}
Let  $W$ be a linear subspace of $V$ and $u\in M$, where $M$ is
defined by (\ref{eqEmb1}). Then the set $M_{u,W}=u+M|_{W}$ is an MDS
code.
\end{proposition}
Proof. By the definitions $M$ is a linear code and $M|_{W}$ is a
linear subspace of $M$. Then $M$ contains an affine subspace
$M_{u,W}=u+M|_{W}$.
 The
code distances  of $M|_{W}$ and $M_{u,W}$ are the same as for $M$.
Firstly, we estimate the cardinality of $M|_{W}$. Given an
$s$-dimensional subspace $W$, there exists a non-degenerate matrix
$B_W$ such that $W=\{xB_W : x\in (GF(q))^s\}$. The rows of $B_W$ are
base of $W$. Let $Z$ be a $d\times s$ matrix over $GF(q)$ such that
$A_MZ=\bar 0$. Then  the matrix $ZB_W$ belongs to $M\cap W=M|_{W}$.
By the McNeish's theorem, the number of such matrices $Z$ is equal
to $(q^s)^{d-t}$. The cardinality of $M_{u,W}$ is the same as one of
$M|_{W}$.
 Moreover, alphabets of any coordinate of elements of $M_{u,W}$ have
the same cardinalities as ones of $M|_{W}$. Then $M_{u,W}$ is an MDS
subcode of $M$ by definition. $\blacksquare$

\section{Proof of Theorem 1}

By Proposition \ref{proNMb}, we can find a prime power $p\geq q$
such that there exists a linear MDS code $M$ over $GF(p)$ of length
$d$ and code distance $t+1$.  Suppose that $A=(I|A')$ is  the check
matrix of this code over $GF(p)$.   Consider elements of $C$. Let
$s_i$ be a number of different symbols in position $i$,
$i=1,\dots,d$. We  use alphabet $V=(GF(p))^n$ where $n=\sum_i s_i$.
Without loss of generality we suppose that $C\subset (GF(p))^n\times
\dots \times (GF(p))^n=V^d$, moreover, we claim that all coordinates
of $d$-tuples from $C$ are unit vectors
$e_i=(0,\dots,0,1,0,\dots,0)\in V$ and the sets of unit vectors in
different positions are disjoint.

 Let $v\in (GF(p))^d$
and $\overline{w}\in V^d$. Define $v\cdot \overline{w}= \sum_i
v_iw_i \in V$. Consider $\overline{w}=(e_1,\dots,e_d)^{T}\in C$. For
each  row $a^i$, $i=1,\dots,t$, of the matrix $A=(I|A')$ we define
vector $g^i=a^i\cdot \overline{w}$. Let
$W=W(w)=L(g^1,\dots,g^t)\subset V$ be a linear hull of vectors
$g^1,\dots,g^t$. By Proposition \ref{proNMe}, the set $ M=\{x\in V^d
: Ax=\bar 0\} $ is an MDS code. It is easy to see that
$\overline{u}=(-g^1+e_1,\dots,-g^t+e_t,e_{t+1},\dots,e_d)^{T}$
belongs to $M$ because
$$A\overline{u}=-A(g^1,\dots,g^t,0,\dots,0)^{T}+A\overline{w}=
-(g^1,\dots,g^t)^{T}+(a^1\overline{w},\dots,a^t\overline{w})^{T}=
\overline{0}.$$ By Proposition \ref{proNMg}, the affine subspace
$\overline{u}+M|_W$ is an MDS subcode of $M$. Moreover,  for each
$i=1,\dots,d$ an alphabet of $i$th coordinate of $\overline{u}+M|_W$
contains $e_i$. By Corollary \ref{proNMf}, there is a switching of
$M$ such that the resulting MDS code $M'$ contains $\overline{w}$.

At last, let us  prove that we can independently make such
switchings for all $\overline{w}\in C$. By Proposition \ref{proNMh},
it is sufficient to prove that subcodes $\overline{u}+M|_W$ are
pairwise disjoint. Let us show that we can reconstruct initial
$\overline{w}$ from $\overline{u}+M|_W$. By definition, every
nonzero linear combination of $g^i$ is a linear combination of rows
of $A$ multiplying by $\overline{w}$. By Corollary \ref{proNMt} and
definition of $e_i$, this linear combination contains $d-t+1$
different vectors $e_i$. If $\overline{x}\in W^d$ contains a nonzero
element in one of the last $d-t$ coordinates then
$\overline{u}+\overline{x}$ contains  at least $d-t$ different
vectors $e_i$ in this coordinate. Since the code distance of $C$ is
equal to $t+1$, we conclude that $\overline{u}+\overline{x}$ belongs
to only one subcode. If $\overline{x}\in W^d$ has zeros  in all last
$d-t$ coordinates then $\overline{u}+\overline{x}=(f^1,\dots,f^t,
e_{t+1},\dots,e_d)^{T}$. Since the code distance of $C$ is equal to
$t+1$  and the Hamming distance between $\overline{u}+\overline{x}$
and $\overline{w}$ is less than $t+1$, we conclude that
$\overline{w}$ is the initial vector for this subcode.

\section{Embedding into Latin hypercube}

\begin{theorem}\label{thEm3}
Let $C\subset Q_q^d$ be a code with code distance $2$. Then there
exist $q',  q\leq q'\leq  q^{d-1}$, and  an $MDS(1,d,q')$ code $M$
such that $C\subset M\subset Q_{q'}^d$.
\end{theorem}

In other words we can embed  a partial Latin $(d-1)$-dimensional
hypercube of order $q$ into  $(d-1)$-dimensional Latin hypercube of
order $q^{d-1}$.

\begin{lemma}\label{lemEmb}
Let $C\subset Q_q^d$ be a code with distance $2$. Suppose that for
each $a\in Q_q$ the retract $C^d_a$ is a subset of an MDS code
$M_a\subset Q_q^{d}$. Then there exists an MDS code $M\subset
Q_{q^2}^{d}$ such that $C\subseteq M$.
\end{lemma}
Proof. Consider an arbitrary MDS code $B\subset Q_q^d$ with code
distance $2$. By Proposition \ref{Cartprod1}, there exists  an MDS
code
$$M'=\{((x_1,y_1),\dots,(x_d,y_d)) : x\in B, y\in
M_a \mbox{ \rm if}\ x_d=a \}$$ constructed as the generalized
Cartesian product. For  $a\in Q_q$ and
$\overline{z}=(z_1,\dots,z_{d-1},a)\in M_a$ we consider the subcube
$E_{\overline{z}}=(Q_q,z_1)\times\cdots\times (Q_q,z_{d-1})\times
H$, where $H=\{(a,a) : a\in Q_q\}$. The set $E_{\overline{z}}\cap
M'$ is an MDS subcode because all  retracts
$$\widehat{(E_{\overline{z}}\cap M')^d_{(a,a)}}=\widehat{B^n_a}\times
{(z_1,\dots,z_{d-1})}$$ are disjoint MDS codes under alphabet
$Q_q\times z_i$ on $i$th coordinate, where $i=1,\dots,d-1$. It is
easy to see that $\widetilde{z}=((0,z_1),\dots,(0,z_{d-1}),(a,a))\in
E_{\overline{z}}$. By Proposition \ref{proNMa} and Corollary
\ref{proNMf}, there exists an MDS subcode $D_{\overline{z}}$ such
that $M''=(M'\setminus (E_{\overline{z}}\cap M'))\cup
D_{\overline{z}}$ is an MDS code and $\widetilde{z}\in M''$. Since
for different $\overline{z}\in C$ subcubes $E_{\overline{z}}$ are
disjoint, all such switchings are independent. By switchings all
$\overline{z}\in C$, we obtain an MDS code $M$ containing
$((0,z_1),\dots,(0,z_{d-1}),(a,a))$ for every
$(z_1,\dots,z_{d-1},a)\in C$.  $\blacksquare$

Proof of Theorem \ref{thEm3}. We use induction on $d$. For $d=2$ a
code $C$ is a partial permutation. It is easy to see that any
partial permutation is embedded into permutation. Suppose that the
theorem is true for $d$. Consider a code $C\subset Q_q^{d+1}$ with
distance $2$. A retract $\widehat{C^{d+1}_a}\subset Q_q^{d}$ is a
code with distance $2$. By the induction hypothesis, there exists an
MDS code $\widehat{M_a} \subset Q_q^{d}$ that includes
$\widehat{C^{d+1}_a}$. By Proposition \ref{proNMq}, for every
$d$-dimensional MDS code $B$  with distance $2$ there exists a
$(d+1)$-dimensional MDS code with distance $2$ which includes $B$ as
a retract. Then for each $a\in Q_q$ we obtain a $(d+1)$-dimensional
MDS code $M_a$ containing $\widehat{M_a} \subset Q_q^{d}$ as a
retract and $C^{d+1}_a$ as the subset. It remains to apply Lemma
\ref{lemEmb} to complete the induction step. $\blacksquare$

For example, consider a partial Latin square of order $3$
$C_3=\begin{array}{|c|c|c|}
\hline {\bf 0}&{\bf 3}&{\bf 6} \\
\hline
{\bf 3}&{ \,}&{\bf 0} \\
 \hline
 \, &{\bf 6}&\, \\
 \hline
\end{array}$ .
By the construction from Theorem \ref{thEm3}, we can embed $C_3$
into a Latin square of order $9$. Consider generalized Cartesian
product of $A=\begin{array}{|c|c|c|}
\hline {a}&{ b}&c \\
\hline
b&{ c}&{ a}\\
 \hline
 c&a&b\\
\hline
\end{array}$\,
and Latin squares $U_a=\begin{array}{|c|c|c|}
\hline {\bf 0}&{ 1}&2 \\
\hline
1&{ 2}&{\bf 0}\\
 \hline
 2&0&1\\
\hline
\end{array}$ ,
$U_c=\begin{array}{|c|c|c|}
\hline 8&{ 7}&{\bf 6} \\
\hline
6&{8}&7 \\
 \hline
7&{\bf 6}&{ 8} \\
\hline
\end{array}$ ,
$U_b=\begin{array}{|c|c|c|}
\hline {4}&{\bf 3}&5 \\
\hline
{\bf 3}&5&4 \\
 \hline
5&4&{ 3}\\
\hline
\end{array}$\,.
\vskip5mm

Elements of the target partial Latin square  are bold. Elements of
the MDS subcode which we choose for switching are  italic. We need
to perform $4$ independent switchings. \vskip5mm

$L_1=\begin{array}{|c|c|c|c|c|c|c|c|c|}
\hline {\bf 0}&{\it 1}&2&4&{\it 3}&{ 5}&8&{ \it 7}&{ 6} \\
\hline
1&{ 2}&{\bf 0}&{3}&5&{4}&{6}&8&{ 7} \\
 \hline
 2&0&1&5&4&3&7&6&8\\
\hline 4&{\it 3}&5&{8}&{\it 7}&{\bf 6}&0&{\it 1}&{ 2} \\
\hline
{ 3}&5&4&6&{8}&7& 1&2&{0} \\
 \hline
5&{ 4}&3&7&{\bf 6}&{ 8}&{ 2}&0&1 \\
\hline 8&{\it 7}&6&0&{\it 1}&{ 2}&{4}&{\textit {\textbf 3}}&5 \\
\hline
6&8&7&1&2&0&{\bf  3}&5&4 \\
 \hline
{ 7}&{ 6}&8& 2&{ 0}&1&5&4&{ 3}\\
\hline
\end{array}$ \hskip10mm $\Rightarrow$ \hskip10mm
$L_2=\begin{array}{|c|c|c|c|c|c|c|c|c|}
\hline {\bf 0}&{\textit {\textbf 3}}&2&4&{\it 7}&{ 5}&8&{\it 1}&{ 6} \\
\hline
1&{ 2}&{\bf 0}&{3}&5&{4}&{6}&8&{ 7} \\
 \hline
 2&0&1&5&4&3&7&6&8\\
\hline 4&{\it 7}&5&{8}&{\it 1}&{\bf 6}&0&{\it 3}&{ 2} \\
\hline
{ 3}&5&4&6&{8}&7& 1&2&{0} \\
 \hline
5&{ 4}&3&7&{\bf 6}&{ 8}&{ 2}&0&1 \\
\hline 8&{\it 1}&6&0&{\it 3}&{ 2}&{4}&{\it 7}&5 \\
\hline
6&8&7&1&2&0&{\bf 3}&5&4 \\
 \hline
{ 7}&{ 6}&8& 2&{ 0}&1&5&4&{ 3}\\
\hline
\end{array}$

\vskip5mm

$L_2=\begin{array}{|c|c|c|c|c|c|c|c|c|}
\hline {\bf 0}&{\bf 3}&2&4&{ 7}&{ 5}&8&{ 1}&{ 6} \\
\hline
{\it 1}&{ 2}&{\bf 0}&{\it 3}&5&{4}&{\it6}&8&{ 7} \\
 \hline
 2&0&1&5&4&3&7&6&8\\
\hline 4&{ 7}&5&{8}&{\it 1}&{\bf 6}&0&{ 3}&{ 2} \\
\hline
{\it 3}&5&4&{\it 6}&{8}&7& {\it 1}&2&{0} \\
 \hline
5&{ 4}&3&7&{\bf 6}&{ 8}&{ 2}&0&1 \\
\hline 8&{1}&6&0&{ 3}&{ 2}&{4}&{ 7}&5 \\
\hline
{\it 6}&8&7&{\it 1}&2&0&{\textit {\textbf 3}}&5&4 \\
 \hline
{ 7}&{ 6}&8& 2&{ 0}&1&5&4&{ 3}\\
\hline
\end{array}$\hskip10mm $\Rightarrow$ \hskip10mm
$L_3=\begin{array}{|c|c|c|c|c|c|c|c|c|}
\hline {\bf 0}&{\bf 3}&2&4&{ 7}&{ 5}&8&{ 1}&{ 6} \\
\hline
{\textit {\textbf 3}}&{ 2}&{\bf 0}&{\it 6}&5&{4}&{\it1}&8&{ 7} \\
 \hline
 2&0&1&5&4&3&7&6&8\\
\hline 4&{ 7}&5&{8}&{ 1}&{\bf 6}&0&{ 3}&{ 2} \\
\hline
{\it 1}&5&4&{\it 3}&{8}&7& {\it 6}&2&{0} \\
 \hline
5&{ 4}&3&7&{\bf 6}&{ 8}&{ 2}&0&1 \\
\hline 8&{1}&6&0&{ 3}&{ 2}&{4}&{ 7}&5 \\
\hline
{\it 6}&8&7&{\it 1}&2&0&{\it 3}&5&4 \\
 \hline
{ 7}&{ 6}&8& 2&{ 0}&1&5&4&{ 3}\\
\hline
\end{array}$

\vskip5mm $L_3=\begin{array}{|c|c|c|c|c|c|c|c|c|}
\hline {\bf 0}&{\bf 3}&{\it 2}&4&{ 7}&{\it 5}& 8&{ 1}&{\it 6} \\
\hline
{\bf 3}&{ 2}&{\bf 0}&{ 6}&5&{4}&{1}&8&{ 7} \\
 \hline
 2&0&1&5&4&3&7&6&8\\
\hline 4&{ 7}&{\it 5}&{8}&{ 1}&{\textit {\textbf 6}}&0&{ 3}&{\it 2} \\
\hline
{ 1}&5&4&{ 3}&{8}&7& { 6}&2&{0} \\
 \hline
5&{ 4}&3&7&{\bf 6}&{ 8}&{ 2}&0&1 \\
\hline 8&{1}&{\it 6}&0&{ 3}&{\it 2}&{4}&{ 7}&{\it 5} \\
\hline
{ 6}&8&7&{1}&2&0&{ 3}&5&4 \\
 \hline
{ 7}&{ 6}&8& 2&{ 0}&1&5&4&{ 3}\\
\hline
\end{array}$
\hskip10mm $\Rightarrow$ \hskip10mm
 $L_4=\begin{array}{|c|c|c|c|c|c|c|c|c|}
\hline {\bf 0}&{\bf 3}&{\textit {\textbf 6}}&4&{ 7}&{\it 2}& 8&{ 1}&{\it 5} \\
\hline
{\bf 3}&{ 2}&{\bf 0}&{ 6}&5&{4}&{1}&8&{ 7} \\
 \hline
 2&0&1&5&4&3&7&6&8\\
\hline 4&{ 7}&{\it 2}&{8}&{ 1}&{\it 5}&0&{ 3}&{\it 6} \\
\hline
{ 1}&5&4&{ 3}&{8}&7& { 6}&2&{0} \\
 \hline
5&{ 4}&3&7&{\bf 6}&{ 8}&{ 2}&0&1 \\
\hline 8&{1}&{\it 5}&0&{ 3}&{\it 6}&{4}&{ 7}&{\it 2} \\
\hline
{ 6}&8&7&{1}&2&0&{ 3}&5&4 \\
 \hline
{ 7}&{ 6}&8& 2&{ 0}&1&5&4&{ 3}\\
\hline
\end{array}$

\vskip5mm $L_4=\begin{array}{|c|c|c|c|c|c|c|c|c|}
\hline {\bf 0}&{\bf 3}&{\bf 6}&4&{ 7}&{ 2}& 8&{ 1}&{ 5} \\
\hline
{\bf 3}&{ 2}&{\bf 0}&{ 6}&5&{4}&{1}&8&{ 7} \\
 \hline
 2&{\it 0}&1&5&{\it 4}&3&7&{\it 6}&8\\
\hline 4&{ 7}&{ 2}&{8}&{ 1}&{ 5}&0&{ 3}&{ 6} \\
\hline
{ 1}&5&4&{ 3}&{8}&7& { 6}&2&{0} \\
 \hline
5&{\it 4}&3&7&{\textit {\textbf 6}}&{ 8}&{ 2}&{\it 0}&1 \\
\hline 8&{1}&{ 5}&0&{ 3}&{ 6}&{4}&{ 7}&{ 2} \\
\hline
{ 6}&8&7&{1}&2&0&{ 3}&5&4 \\
 \hline
{ 7}&{\it 6}&8& 2&{\it 0}&1&5&{\it 4}&{ 3}\\
\hline
\end{array}$
\hskip10mm $\Rightarrow$ \hskip10mm
$L_5=\begin{array}{|c|c|c|c|c|c|c|c|c|}
\hline {\bf 0}&{\bf 3}&{\bf 6}&4&{ 7}&{ 2}& 8&{ 1}&{ 5} \\
\hline
{\bf 3}&{ 2}&{\bf 0}&{ 6}&5&{4}&{1}&8&{ 7} \\
 \hline
 2&{\textit {\textbf 6}}&1&5&{\it 0}&3&7&{\it 4}&8\\
\hline 4&{ 7}&{ 2}&{8}&{ 1}&{ 5}&0&{ 3}&{ 6} \\
\hline
{ 1}&5&4&{ 3}&{8}&7& { 6}&2&{0} \\
 \hline
5&{\it 4}&3&7&{\it 6}&{ 8}&{ 2}&{\it 0}&1 \\
\hline 8&{1}&{ 5}&0&{ 3}&{ 6}&{4}&{ 7}&{ 2} \\
\hline
{ 6}&8&7&{1}&2&0&{ 3}&5&4 \\
 \hline
{ 7}&{\it 0}&8& 2&{\it 4}&1&5&{\it 6}&{ 3}\\
\hline
\end{array}$

\end{document}